\documentclass[conference]{IEEEtran}
\usepackage{ifpdf}
\usepackage{cite}
\ifCLASSINFOpdf
  \usepackage[pdftex]{graphicx}
  \graphicspath{{C:\Users\diat\Desktop\Conference 2016\IEEE_conf/pdf/}}

  \DeclareGraphicsExtensions{.pdf,.jpeg,.png}
\else
   \usepackage[dvips]{graphicx}
\fi
\usepackage{amsmath}
\usepackage{graphicx}
\usepackage{xcolor}
\usepackage{algorithm,algcompatible,amsmath}
\usepackage{amsmath}
\usepackage{amssymb}
\usepackage{subfig}
\usepackage{multicol,lipsum}

\usepackage{array}

\usepackage{multicol, blindtext} 
\usepackage{fixltx2e}
\begin{document}
\title{Modified Recursive Cholesky (Rchol) Algorithm  \newline
	{\huge An Explicit Estimation and Pseudo-inverse of Correlation Matrices}}

\author{\IEEEauthorblockN{Vanita Pawar}
\IEEEauthorblockA{
 vanita\_phdee12@diat.ac.in}
\and
\IEEEauthorblockN{Krishna Naik Karamtot}
\IEEEauthorblockA{
 krishnanaik@diat.ac.in}
}


%Telephone: (800) 555--1212, Fax: (888) 555--1212}


\maketitle

\begin{abstract}
The Cholesky decomposition plays an important role in finding the inverse of the correlation matrices. As it is a fast and numerically stable for linear system solving, inversion, and factorization compared to singular valued decomposition (SVD), QR factorization and LU decomposition. As different methods exist to find the Cholesky decomposition of a given matrix, this paper presents the comparative study of a proposed RChol algorithm with the conventional methods. The RChol algorithm is an explicit way to estimate the modified Cholesky factors of a dynamic correlation matrix.

\end{abstract}

\vspace{.5in}

%\begin{IEEEkeywords}
%	Modified Cholescky decomposition, BEM, Linear Prediction, Multipath
%\end{IEEEkeywords}

Cholesky decomposition is a fast and numerically stable for linear system solving, inversion, and factorization compared to singular valued decomposition (SVD), QR factorization and LU decomposition $[1]$.  The wireless communication system is highly dependent on matrix inversion of the correlation matrix. Such system consists of a huge matrix inversion. An outdoor wireless communication has a time-varying channel which changes dynamically for mobile user. In case of narrowband channel, the channel is considered constant for a symbol duration, whereas for broadband, it is changing within a symbol period. Such time-varying channel forms the special structure of channel matrix and correlation matrix. To exploit such special structure, a novel modified recessive Cholesky (RChol) algorithm is introduced in $[2]$. Our proposed (RChol) algorithm is a computational efficient algorithm to compute the modified Cholesky factors of known as well as an unknown covariance matrix. 

In this paper, we present the comparative study of conventional Cholesky algorithm and the RChol algorithm to manifest the importance of the proposed algorithm in a highly dynamic wireless communication.
%{\color{red}{Insert antenna Diversity}}
\section{System Model}
In wireless communication system, number of transmit and or received antennas are used to improve   the diversity of the system. The channel $‘h’$ between transmitter and receiver has the different form and depends on the number of antennas used at the transmitter and the receiver side. The channel for Single-input-single-output (SISO) as $h=\{h^{0}_{n}, h^{1}_{n}, \hdots h^{L-1}_{n}\}$, for	 Single-input-multiple-output (SIMO) as ${\bf h}=\{{\bf h}^{0}_{n}, {\bf h}^{1}_{n}, \hdots {\bf h}^{L-1}_{n}\}$ and for	Multiple-input-multiple-output (MIMO): ${\bf H}=\{{\bf H}^{0}_{n}, {\bf H}^{1}_{n}, \hdots {\bf H}^{L-1}_{n}\}$.
%\begin{figure}[h]
%	\centering
%	{\includegraphics[width=60mm]{Antenna_diversity.png}} 
%	\caption{Space Diversity  (a) SISO, (b) SIMO, (c) $SuMIMO$, (d) $MuMIMO$ } \label{fig:Ant_Div}
%\end{figure} 

Let {\bf y}(n) be received signal  with the number of transmit antennas $'K=1'$, multipath $‘L-1’$ and channel noise $‘v’$, represented as

\begin{align}
{\bf y}(n):=\sum^{K}_{k=1} \sum^{L-1}_{l=0}{\bf h}^{k}(n;l)s^{k}(n-l)+{\bf v}(n) ,\hspace{5mm} n=0,1,..T-1
\label{eq:Mi_y}
\end{align}

Let ${\bf y}_{N}(n)$ be the received vector by stacking $'N'$ successive received vectors.
Where  ${\bf y}_{N}(n)=[{\bf y}(n), {\bf y}(n-1), \hdots {\bf y}(n-N+1) ]^T$ and the transmitted symbol vector is ${\bf s}_{N}=[ s(n),s(n-1),….s(n-N+1) ]^T$.  
%and ${\bf H}_{N}(n)$ as equation ({\ref{eq:H_mat}}).
%
%
%	\begin{align}{
%		\left [ \begin{array}{cccccccc}
%		{\bf H}^{0}_{n}&{\bf H}^{1}_{n}&\cdots&{\bf H}^{L-1}_{n}&{\bf 0}&\cdots&{\bf 0} \\
%		{\bf 0 }& {\bf H}^{0}_{n-1} & {\bf H}^{1}_{n-1} & \cdots & {\bf H}^{L-1}_{n-1}& \cdots & {\bf 0}\\
%		\vdots & \vdots & \vdots & \vdots & \vdots & \ddots & \vdots  \\
%		{\bf 0} & {\bf 0} & \cdots &{\bf H}^{0}_{n-N+1} & {\bf H}^{1}_{n-N+1} & \cdots & {\bf H}^{L-1}_{n-N+1}  \\
%		\end{array} \right ]}
%	\label{eq:H_mat}
%	\end{align}
%	
%${\bf H}_{N}(n-1)$ as equation ({\ref{eq:H_mat1}}).
%
%	\begin{align}{
%	\left [ \begin{array}{cccccccc}
%	{\bf H}^{0}_{n-1}&{\bf H}^{1}_{n-1}&\cdots&{\bf H}^{L-1}_{n-1}&{\bf 0}&\cdots&{\bf 0} \\
%	{\bf 0 }& {\bf H}^{0}_{n-2} & {\bf H}^{1}_{n-2} & \cdots & {\bf H}^{L-1}_{n-2}& \cdots & {\bf 0}\\
%	\vdots & \vdots & \vdots & \vdots & \vdots & \ddots & \vdots  \\
%	{\bf 0} & {\bf 0} & \cdots &{\bf H}^{0}_{n-N} & {\bf H}^{1}_{n-N} & \cdots & {\bf H}^{L-1}_{n-N}  \\
%	\end{array} \right ]}
%\label{eq:H_mat1}
%\end{align}
%\begin{figure*}[]	
%	%and  ${\bf H}_{N}(n-1)$ as ,
%	\begin{equation}{
%		\small
%		{\bf H}_{N}(n-1)=
%		\left [ \begin{array}{ccccccc}
%		{\bf H}(n-1;0) & {\bf H}(n-1;1) & \cdots & {\bf H}(n-1;L-1) & {\bf 0} & \cdots & {\bf 0} \\
%		{\bf 0 }& {\bf H}(n-2;0) & {\bf H}(n-2;1) & \cdots & {\bf H}(n-2;L-1)& \cdots & {\bf 0}\\
%		\vdots & \vdots & \vdots & \vdots & \vdots & \ddots & \vdots  \\
%		{\bf 0} & {\bf 0} & \cdots &{\bf H}(n-N;0) & {\bf H}(n-N;1) & \cdots & {\bf H}(n-N;L-1)  \\
%		\end{array} \right ]}
%	\label{eq:H_mat_nn}
%	\end{equation}
%\end{figure*}
Then ${\bf y}_{N}(n)$ can be represented in matrix form as ${\bf y}_N(n)={\bf H}_N{\bf s_N(n)+{\bf v}_N(n)}$
%\begin{align}
%{\bf y}_N(n)={\bf H}_N{\bf s_N(n)+{\bf v}_N(n)}
%\label{eq:Y_N}
%\end{align}
and the correlation matrix for  ${\bf y}_N$ can be written as ${\bf R}_N (n)=E[{\bf y}_N (n) {\bf y}_N^H (n)]$. Let ${\bf r}^{n}_{00}=E[{\bf y}(n){\bf y}^H(n)]$ and ${\bf r}^{n}_{ij}=E[{\bf y}(n-i){\bf y}^H (n-j)]$ then the correlation matrix ${\bf R}_N (n)$ and ${\bf R}_N(n-1)$  at time instant $'n'$ and $'n-1'$ can be represented as equation (\ref{eq:R1}) and equation (\ref{eq:R2}) respectively. 

\begin{equation}{
	\tiny
	\left[ \begin{array}{ccccc}
	{\bf r}^{n}_{00}&{\bf r}^{n}_{01}&\hdots \hdots& {\bf r}^{n}_{0(N-2)}&{\bf r}^{n}_{0(N-1)} \\
	{\bf r}^{n}_{10} & {\color{cyan}{\bf r}^{n}_{11}}&  {\color{cyan}\hdots \hdots}&  {\color{cyan} {\bf r}^{n}_{1(N-2)}}& {\color{cyan}{\bf r}^{n}_{1(N-1)}}\\
	\vdots&  {\color{cyan}\vdots}& {\color{cyan}\ddots}& {\color{cyan}\vdots}& {\color{cyan}\vdots}\\
	%	{\bf r}^{n}_{(N-2)0}& {\color{cyan}{\bf r}^{n}_{(N-2)1}}& {\color{cyan}\hdots \hdots}&  {\color{cyan}{\bf r}^{n}_{(N-2)(N-2)}}& {\color{cyan}{\bf r}^{n}_{(N-2)(N-1)}} \\
	{\bf r}^{n}_{(N-1)0}& {\color{cyan}{\bf r}^{n}_{(N-1)1}}& {\color{cyan}\hdots \hdots}&  {\color{cyan}{\bf r}^{n}_{(N-1)(N-2)}}& {\color{cyan}{\bf r}^{n}_{(N-1)(N-1)}}\\
	\end{array} \right]
	\label{eq:R1}
}\end{equation}

\begin{equation}{
	\tiny
	\left[ \begin{array}{ccccc}
	{\color{cyan}{\bf r}^{n}_{11}}& {\color{cyan}{\bf r}^{n}_{12}}& {\color{cyan}\hdots \hdots}&  {\color{cyan}{\bf r}^{n}_{1(N-1)}}& {\bf r}^{n}_{1N} \\
	%	{\color{cyan}{\bf r}^{n}_{21}}& {\color{cyan}{\bf r}^{n}_{22}}& {\color{cyan}\hdots \hdots}&  {\color{cyan}{\bf r}^{n}_{2(N-1)}}&{\bf r}^{n}_{2N}\\
	{\color{cyan}\vdots}&{\color{cyan}\vdots}&{\color{cyan}\ddots}&{\color{cyan}\vdots}&\vdots\\
	{\color{cyan}{\bf r}^{n}_{(N-1)1}}& {\color{cyan}{\bf r}^{n}_{(N-1)2}}& {\color{cyan}\hdots \hdots}&  {\color{cyan}{\bf r}^{n}_{(N-1)(N-1)}}&{\bf r}^{n}_{(N-1)N} \\
	{\bf r}^{n}_{N1}&{\bf r}^{n}_{N2}&\hdots \hdots& {\bf r}^{n}_{N(N-1)}&{\bf r}^{n}_{NN}\\
	\end{array} \right]
	\label{eq:R2}
}\end{equation}

%and  R.H.S can be as in equation (\ref{R_E}).
%  
%%\begin{equation}{
%%	\small
%%	{\bf R}_N (n)=
%%	\left[ \begin{array}{ccccc}
%%	E[{\bf y}(n){\bf y}^H(n)]  &E[{\bf y}(n){\bf y}^H(n-1)] & \hdots& E[{\bf y}(n){\bf y}^H(n-N)]\\
%%	E[{\bf y}(n-1){\bf y}^H(n)]  &E[{\bf y}(n-1){\bf y}^H(n-1)] & \hdots &E[{\bf y}(n-1){\bf y}^H(n-N)]\\
%%	\vdots& \vdots &\ddots & \vdots\\
%%	E[{\bf y}(n-N){\bf y}^H(n)]  &E[{\bf y}(n-N){\bf y}^H(n-1)] & \hdots &E[{\bf y}(n-N){\bf y}^H(n-N)]\\
%%	\end{array} \right]
%%	\label{R_E}
%%}\end{equation}
%\begin{equation}{
%	\small
%%	{\bf R}_N (n)=
%	\left[ \begin{array}{ccc}
%	E[{\bf y}(n){\bf y}^H(n)]   & \hdots& E[{\bf y}(n){\bf y}^H(n-N)]\\
%%	E[{\bf y}(n-1){\bf y}^H(n)]  &E[{\bf y}(n-1){\bf y}^H(n-1)] & \hdots &E[{\bf y}(n-1){\bf y}^H(n-N)]\\
%	\vdots &\ddots & \vdots\\
%	E[{\bf y}(n-N){\bf y}^H(n)] & \hdots &E[{\bf y}(n-N){\bf y}^H(n-N)]\\
%	\end{array} \right]
%	\label{R_E}
%}\end{equation}
%\begin{align}
%{\bf R}_N (n)=[{\bf R}^{c_1} (n),{\bf R}^{c_2} (n),\hdots {\bf R}^{c_N} (n)]
%\label{R_clm}
%\end{align}

\section{Cholesky Decomposition}
The correlation matrix is complex matrix and the pseudo-inverse of {\bf R} can be computed from Cholesky factors, such that if lower triangular matrix $‘{\bf L}’$ is Cholesky factors of the correlation matrix $‘{\bf R}’$ and can be represented as  ${\bf R}={\bf LL}^H$  then pseudo-inverse of {\bf R} can be computed as ${\hat{ \bf R}}^{\dagger}:={\bf L}^{-H}{\bf L}^{-1}$.
The section below details the conventional Cholesky algorithms and the RChol algorithm.

\subsection{Cholesky Decomposition (Gaxpy version)}
The Cholesky Decomposition $[3]$ factorizes a complex (or real-valued) positive-deﬁnite Hermitian symmetric matrix into a product of a lower triangular matrix and its Hermitian transpose. ${\bf R}={\bf LL}^H$  where, {\bf L} is a lower triangular matrix and ${\bf L}^H$ is Hermitian of {\bf L}. The matrix {\bf R} must be a positive definite and this method needs square root operation.
\subsubsection{Algorithm steps}
\begin{enumerate}
	\item  Compute {\bf R} at each time instant n
	\item Find the square root of diagonal element of {\bf R}
	\item Modify each column of {\bf R}
	\item Equate lower triangular part of {\bf R} to {\bf L}
	\item Repeat steps  $(1)$ to $(4)$ for each time instant
\end{enumerate}

\begin{algorithm}[h]
	\tiny {
		\caption{\small \bf {Cholesky Decomposition} ${\bf R}={\bf LL}^H$
			%						\cite{Golub2012},  \cite{RHunger}\hspace{3mm} ${\bf R}={\bf LL}^H $
		}
		\begin{algorithmic}
			\STATE \textbf{Initialization:}
			%\newline   for $k=1$ \[{\bf H}_{1}(n)=[{\bf r}_{00},{\bf  r}_{10},.....{\bf r}_{(N-1)0}]^{T}\]
			%\[{\bf {\tilde H}}_{1}(n)=[{\bf 0},{\bf r}_{10},.....{\bf r}_{(N-1)0}]^{T}\]
			%					\newline  
			\[[{\bf R}]_{1:N,1}=\frac{[{\bf R}]_{1:N,1}}{\color{cyan}\sqrt{[{\bf R}]_{1,1}}}\]
			\STATE \textbf{Order Updates on ${ \bf R}^{'}s$:} 
			%					\newline
			\colorbox{yellow}{\text {for $k=2$ to $N$ }}
			\[[{\bf R}]_{{k:N},k}=[{\bf R}]_{{k:N},k}-[{\bf R}]_{{k:N},{1:k-1}}[{\bf R}]^{H}_{{k},{1:k-1}}\]
			\[[{\bf R}]_{k:N,k}=\frac{[{\bf R}]_{k:N,k}}{\color{cyan}\sqrt{[{\bf R}]_{k,k}}}\]
			end
			\[{\bf L}=\text {tril}({\bf R})\]
		\end{algorithmic}
	}
\end{algorithm}

\subsection{Modified Cholesky Algorithm ${\bf R}={\bf LDL}^H$}
To avoid square root operation, a modified Cholesky algorithm $[3]$ is used, which avoids square root operation by introducing a diagonal matrix D in between Cholesky factors. The modified Cholesky algorithm does not require {\bf R} to be a positive definite matrix but it's determinant must be nonzero.  {\bf R} may be rank deficient to a certain degree i.e. D may contain negative main diagonal entries if {\bf R} is not positive semi-definite.
\subsubsection{Algorithm steps}
\begin{enumerate}
	\item Compute {\bf R} at each time instant n
	\item Modify each column of {\bf R}
	\item Equate the strictly lower part of matrix {\bf R} to ${\bf L}_1$ with ones on the main diagonal
	\item	Equate main diagonal of {\bf R} with the main diagonal of {\bf D}
	\item Repeat step $(1)$ to $(4)$ for each time instant.
\end{enumerate}

\begin{algorithm}[h]
	\tiny {
		\caption{\small {\bf Modified Cholesky Decomposition} ${\bf R}={\bf LDL}^H$ 
			%						 \cite{Golub2012}, \cite{RHunger} ${\bf R}={\bf LDL}^H$ 
		}
		\begin{algorithmic}
			\STATE \textbf{Initialization:}
			%\newline   for $k=1$ \[{\bf H}_{1}(n)=[{\bf r}_{00},{\bf  r}_{10},.....{\bf r}_{(N-1)0}]^{T}\]
			%\[{\bf {\tilde H}}_{1}(n)=[{\bf 0},{\bf r}_{10},.....{\bf r}_{(N-1)0}]^{T}\]
			%					\newline  
			\[[{\bf R}]_{2:N,1}=\frac{[{\bf R}]_{2:N,1}}{[{\bf R}]_{1,1}}\]
			\STATE \textbf{Order Updates on ${ \bf R}^{'}s$:} 
			%					\newline
			\colorbox{yellow}{\text {for $k=2$ to $N$ }}\\
			\colorbox{green}{for i=1:k-1}\\
			\hspace{9mm}	$[{\textit{\bf v}}]_i=[{\bf R}]_{1,k}$ \hspace{3mm} , \text {if $i=1$} \\
			\hspace{9mm}   $[{\bf R}]_{i,i}[{\bf R}]^{*}_{n,i}$ ,  \text {if $i\neq1$} \\
			\colorbox{green}{end}
			\[[{\textit{\bf v}}]_k=[{\bf R}]_{k,k}-[{\bf R}]_{k,{1:k-1}}[{\textit{\bf v}}]_{1:k-1}\]
			\[[{\bf R}]_{k,k}=[{\textit{\bf v}}]_k\]
			\[[{\bf R}]_{{k+1:N},k}=\frac{[{\bf R}]_{{k+1:N},k}-[{\bf R}]_{{k+1:N},{1:k-1}}[{\textit{\bf v}}]_{1:k-1}}{[{\textit{\bf v}}]_k}\]
			end\\
			\[{\bf D} = diag(daig({\bf R}))\]
			\[{\bf L}=\text {tril}({\bf R})\]
		\end{algorithmic}
	}
\end{algorithm}

\subsection{Recursive  Cholesky Algorithm (The Shcur Algorithm) ${\bf R}_{Schur}:={\bf HH}^H $}

The Schur algorithm recursively compute the columns of the lower triangular matrix {\bf H} form matrix {\bf R}. It is shown in $[4]$ that Levinson recursion may be used to derive the Lattice recursion for computing QR factors of data matrices and Lattice recursion can be used to derive the Schur recursion for computing Cholesky factors of a Toeplitz correlation matrix. The detail algorithm is given in algorithm $3$. The Schur algorithm like previously mentioned algorithm computes all $‘N’$ inner product to compute matrix {\bf R} for initialization.
%\begin{align}
%{\bf R}_{Schur}:={\bf HH}^H 
%\label{eq:RSChur}
%\end{align}
\subsubsection{Algorithm steps}
\begin{enumerate}
	\item Compute {\bf R} at each time instant n
	\item Initialize first column of {\bf R} to the first column of Cholesky factor {\bf H} 
	\item Compute rest column recursively from columns of {\bf R}
	\item repeat step $(1)$ to $(3)$ for each time instant
\end{enumerate}

\begin{algorithm}[h]
	\tiny{
		\caption{\small {\bf Schur Algorithm} ${\bf R}={\bf LL}^H$
			%						{\cite{scarf} ${\bf R}={\bf LDL}^H$ }
		}
		\begin{algorithmic}
			\STATE \textbf{Initialization:}
			%\newline   for $k=1$ \[{\bf H}_{1}(n)=[{\bf r}_{00},{\bf  r}_{10},.....{\bf r}_{(N-1)0}]^{T}\]
			%\[{\bf {\tilde H}}_{1}(n)=[{\bf 0},{\bf r}_{10},.....{\bf r}_{(N-1)0}]^{T}\]
			\newline   for $k=1$ 
			\[{ H}_{1}(n)=[{ r}^n_{00},{  r}^n_{10},.....{ r}^n_{(N-1)0}]^{T}\]
			\[{{\tilde H}}_{1}(n)=[{ 0},{ r}^n_{10},.....{ r}^n_{(N-1)0}]^{T}\]
			\STATE \textbf{Order Updates on ${ H}^{'}s$:} 
			%					\newline
			\colorbox{yellow}{\text {	for $k=2$ to $N$ }}
			%\[ {\bf \sigma}_{n+1}{\bf H}_{n+1}= {\bf \tilde{k}}_{ref_{k}}[{\bf {k}}_{ref_{k}}({\bf \sigma}_n{\bf \tilde H}_n) +{\bf Z}_M({\bf \sigma}_n{\bf H}_n ) ] \]
			%\[ {\bf \sigma}_{n+1}{\bf\tilde H}_{n+1}= {\bf \tilde{k}}_{ref_{k}}[({\bf \sigma}_n{\bf \tilde H}_n)  +{\bf {k}}_{ref_{k}}{\bf Z}_M({\sigma}_n{\bf H}_n ) ] \]
			\[ { \sigma}_{k}{ H}_{k}= { \tilde{k}}_{ref_{k}}[{ {k}}_{ref_{k}}({ \sigma}_{k-1}{ \tilde H}_{k-1}) +{ Z}_M({ \sigma}_{k-1}{ H}_{k-1} ) ] \]
			\[ { \sigma}_{k}{\tilde H}_{k}= { \tilde{k}}_{ref_{k}}[({ \sigma}_n{ \tilde H}_k)  +{ {k}}_{ref_{k}}{ Z}_M({\sigma}_k{ H}_k ) ] \]
			\STATE \textbf{Scaling Factors:} 
			\[{ {k}}_{ref_{k}}= -{\frac{({ \sigma}_k{\tilde{ H}}_k)_{k}}{ ({ \sigma}_{k}{ H}_{k})_{k-1}}} \]
			\[{ \tilde{k}}_{ref_{k}}=  \frac{({\sigma}_k{ H}_k)_{k-1}}{ ({ \sigma}_{k+1}{ H}_{k+1})_{k}} \]
			%${\bf \tilde{k}_{ref}}_{_k}= \frac{{\bf k_{ref}}(n)^{\bf '}{\bf D}_0{(n-1)}}{\tilde{\bf D}_0(n)}$
		\end{algorithmic}
		$*$ Note: Here notation is followed same as in {\cite{4}}  and H represents vector
	}
\end{algorithm}

\subsection{ The RChol Algorithm  ${\hat {\bf L}}{\hat{\bf D}}{\hat{\bf L}}^H:={\bf R}$}

It is clear from above equation (\ref{eq:R1}) and equation (\ref{eq:R2})  that ${\bf R}_N(n)$ can be represented from submatrix of ${\bf R}_N(n-1)$.
To utilize such special structure of correlation matrices, we propose a modified recursive Cholesky algorithm to compute the Cholesky factors recursively. This algorithm is modification of Schur algorithm mentioned above. The more general approach consists of using the Schur algorithm to induce recursion for columns of dynamic {\bf L}. This algorithm does not need N inner products to compute the correlation matrix {\bf R}. The Cholesky factors are computed explicitly such that 
%\begin{align}
%{\hat {\bf L}}{\hat{\bf D}}{\hat{\bf L}}^H:={\bf R}
%\label{eq:}
%\end{align}
Let ${\bf L}_1={\bf LD}^{1/2}$then pseudo-inverse can be computed as ${\hat {\bf R}}^\dagger ={\bf L}_1^{-H}{\bf L}_1^{-1}$
%\begin{align}
%{\hat {\bf R}}^\dagger ={\bf L}_1^{-H}{\bf L}_1^{-1}
%\end{align}
\subsubsection{Algorithm steps}
\begin{enumerate}
	\item Initialize first the first column of Cholesky factor {\bf A} as ${\bf A}_1$
	\item Compute second column recursively from ${\bf A}_1(n)$ and ${\bf A}_1(n-1)$
	\item Substitute sub-matrix ${\bf A}_{2:N-1,2:N-1}(n-1)$ to ${\bf A}_{3: N,3:N}(n)$
	\item Repeat step $(1)$ to $(3)$ for each time instant
\end{enumerate}

In the Schur algorithm, columns of Cholesky factors at time instant $‘n’$ are computed recursively from the correlation matrix at that instant. Whereas in the RChol algorithm first two columns of Cholesky factors at time instant $‘n’$ is computed recursively from previous Cholesky factor and submatrix of that Cholesky factors are updated recursively from previous Cholesky factor i.e. at time instant $‘n-1’$. Conventional Cholesky algorithm mentioned here are introduced for normal matrices whereas proposed matrix is well suited for block matrices and simulations are shown for that only.
%%%%%%%%%%%%%%%%%%%%%%%%%%%%%%%%%%%%%%%%%%%%%%%%%%%%%%%%%%%%%%%
%\begin{minipage}{\textwidth}
%	\begin{minipage}[t]{0.44\textwidth}

%	\end{minipage}
%	\hfill
%	\begin{minipage}[b]{0.48\textwidth}

%	\end{minipage}
%\end{minipage}
%\newline
%\begin{minipage}{\textwidth}
%	\begin{minipage}[b]{0.45\textwidth}

%	\end{minipage}
%	\hfill
%	\begin{minipage}[b]{0.45\textwidth}
\begin{algorithm}[h]
	\tiny {
		\caption{\small{\bf Recursive Cholesky Update : RChol}   ${\bf R}={\bf LDL}^H$ }
		%				\cite{Vanita}
		\begin{algorithmic}
			\STATE \textbf{Initialization:}
			\newline   for $k=1$ , \hspace{3mm} ${\bf D}_1(n)={\bf r}^n_{00}	$
			\[{\bf A}_{1}(n)=[{\bf r}^n_{00},{\bf  r}^n_{10},\hdots {\bf r}^n_{(N-1)0}]^{T}\]
			\[{\bf {\tilde A}}_{1}(n)=[{\bf 0},{\bf r}^n_{10},\hdots{\bf r}^n_{(N-1)0}]^{T}\]
			%\[{\bf D}_1(n)={\bf r}_{00}	\]
			\STATE \textbf{Order Updates on ${\bf A}^{'}s$:} 
			%					\newline    
			\colorbox{yellow}{\text { for $k=2$ }}
			\[	{\bf A}_{k}(n):={\bf Z}_M {\bf A}_{k-1}(n-1)-\tilde{\bf A}_{k-1}(n){\bf\tilde{ k}_{ref}}(n)\]
			\[ {\bf D}_{k}(n)={\bf D}_{k}(n-1)[ {\bf I}_M-{\bf k}_{ref}(n){\bf\tilde{ k}}_{ref}(n)]\]
			\colorbox{yellow}{\text {for $k>2$ }}, \hspace{3mm}    $\tilde{\bf A}_{k-1}(n)=\bf 0$  
			\[{\bf A}_k(n)={\bf Z_M A}_k(n-1) \]
			\[ {\bf D}_k(n)={\bf Z_M D}_k(n-1) \]
			\STATE \textbf{Scaling Factors:} 
			${\bf k}_{ref}(n)=\frac{\tilde{\bf A}_1(n)_{(2,:)}}{{\bf A}_1{(n-1)}_{(1,:)}} $\\
			${\bf \tilde{k}}_{ref}(n)= \frac{{\bf k}_{ref}(n)^{\bf '}{\bf D}_1{(n-1)}}{\tilde{\bf D}_1(n)}$
		\end{algorithmic}
	}
\end{algorithm}
%	\end{minipage}
%\end{minipage}
\section{Simulation results}
%\begin{figure}[th]
%	\centering
%	\subfloat[Proposed Algorithm (Difference)]{\includegraphics[width=0.4\columnwidth]{BlindB1realDiffn15db}} \hfil
%	\subfloat[Schur Algorithm (Difference)]{\includegraphics[width=0.4\columnwidth]{Schur_Blind_Diff.pdf}} \\
%	\vspace*{10pt}
%	\subfloat[Proposed Algorithm (Ratio)]{\includegraphics[width=0.4\columnwidth]{BlindB1realRatio15db}}  \hfil
%	\subfloat[Schur Algorithm (Ratio)]{\includegraphics[width=0.4\columnwidth]{Schur_Blind_Ratio}} 
%	\vspace*{8pt}
%	\caption{ Comparisons of RChol algorithm Vs Schur Algorithm for the ‘unknown’ correlation matrix ‘R’ }
%	\label{fig:},
%\end{figure}

To compare proposed the RChol algorithm with Schur algorithm, we compared the result of both the algorithm with theoretical results.  Fig. $1$. Show the ratio and difference of matrices ${\hat{\bf R}}_N$, ${\hat{\bf R}}_{RChol}$  and ${\hat{\bf R}}_{Schur}$ , when the correlation matrix is unknown. That has the application in blind channel and or data estimation. Fig. 1 (a) and (b) shows the maximum error for the RChol algorithm, $[{\hat{\bf R}}_N-{\hat{\bf R}}_{RChol} ]$ is $‘0.6’$ while for the Schur algorithm,  $[{\hat{\bf R}}_N-{\hat{\bf R}}_{Schur} ]$  is $‘4’$ i.e. nearly 6 times the RChol algorithm. In case of ratio Fig. 1 (a) and (b) shows the maximum ratio for the RChol algorithm,  $[{\hat{\bf R}}_N./{\hat{\bf R}}_{RChol}]$ is $‘45’$ while for the Schur algorithm, $[{\hat{\bf R}}_N./{\hat{\bf R}}_{Schur}]$  is $‘1500’$.

\begin{figure}[h]
	\centering
	{\includegraphics[width=80mm]{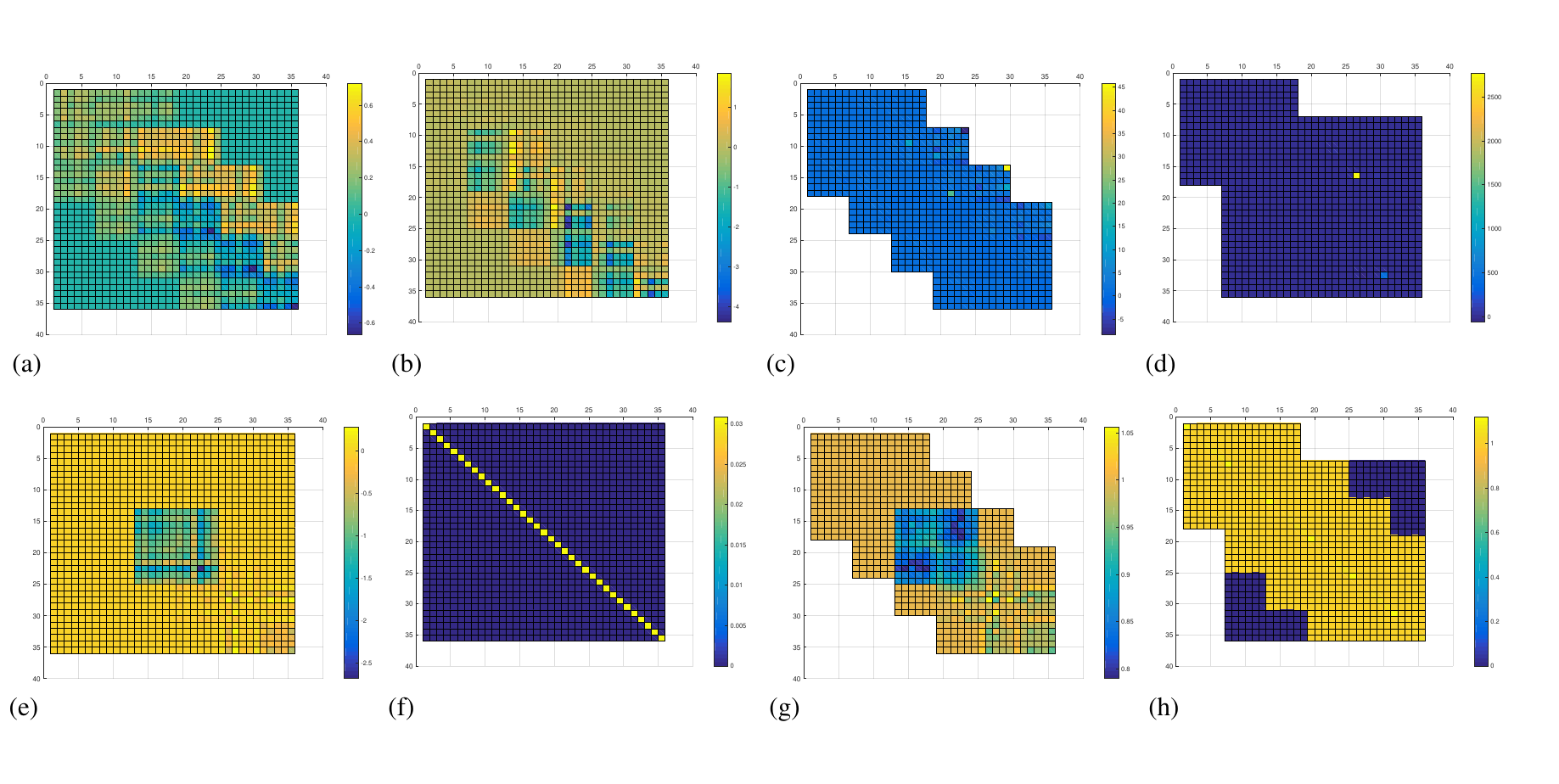}} 
	\caption{ Comparisons of RChol algorithm Vs Schur Algorithm for the $‘unknown’$ and $‘known’$ correlation matrix $‘\bf R’$, (a, e): Proposed Algorithm ( Difference), (b, f): Schur Algorithm( Difference), (c, g):Proposed Algorithm (Ratio), (d, h): Schur Algorithm (Ratio)}
\end{figure} 

Fig. 1 Show the ratio and difference of matrices ${\hat{\bf R}}_N$, ${\hat{\bf R}}_{RChol}$  and ${\hat{\bf R}}_{Schur}$ , when the correlation matrix is known. Fig. 1 (a) and (b) shows that the maximum error for the RChol algorithm,  $[{\hat{\bf R}}_N-{\hat{\bf R}}_{RChol} ]$ is $‘2.5’$ while for the Schur algorithm, $[{\hat{\bf R}}_N-{\hat{\bf R}}_{Schur} ]$  is $‘0.03’$ i.e. nearly 6 times the RChol algorithm. In case of ratio Fig. 1 (e) and (f) shows that the maximum ratio for the RChol algorithm, $[{\hat{\bf R}}_N./{\hat{\bf R}}_{RChol} ]$ is $‘1.15’$ while for the Schur algorithm, $[{\hat{\bf R}}_N./{\hat{\bf R}}_{Schur} ]$ is $‘1’$.

From Fig. 1  it can be concluded that the Schur algorithm is best suited when the correlation matrix is known, but leads to huge error propagation through the column when {\bf R} is unknown and cannot be applied for blind channel estimation. In converse, the RChol algorithm is best suited for blind channel estimation and reduces error propagation through the column.

%\begin{figure}[th]
%%\centerline{\psfig{file=jaaf1.eps,width=4.8cm}}
%\vspace*{8pt}
%\caption{A schematic illustration of dissociative recombination. The
%direct mechanism, 4m$^2_\pi$ is initiated when the molecular ion
%S$_{\rm L}$ captures an electron with kinetic energy.}
%\end{figure}

\section{Conclusion}
Convention  methods  of  Cholesky  factorization  requires  the  correlation  matrix  which  needs  inner 
product.  While  the  recursive  modiﬁed  Cholesky  algorithm  (RChol)  algorithm  is  an  explicit  way  to 
recursively calculating the pseudo-inverse of the matrices without estimating the correlation matrix. It 
requires less number of iteration which avoids error propagation through column updates. The RChol 
algorithm has most of the use in calculating the pseudo-inverse of the of a time-varying matrix which is 
applicable to SIMO/MIMO, CDMA, OFDM, etc. wireless communication systems. 

\vskip5pt

\noindent V. Pawar and K. Naik (\textit{DIAT, Pune, India})
\vskip3pt

\noindent E-mail: vanietaapawar@gmail.com

\end{document}